# Division in Associative $D$-Algebra

Aleks Kleyn

ABSTRACT. From the symmetry between definitions of left and right divisors in associative $D$-algebra $A$, the possibility to define quotient as $A \otimes A$-number follows. In the paper, I considered division and division with remainder. I considered also definition of prime $A$-number.

## Contents

## 1. Preface

Let $D$ be commutative ring. Let $D$-algebra $A$ be associative.

**Definition 1.1.** $A$-number $a$ is left divisor of $A$-number $b$, if there exists $A$-number $c$ such that

$$ac = b \tag{1.1}$$

□

**Definition 1.2.** $A$-number $a$ is right divisor of $A$-number $b$, if there exists $A$-number $c$ such that

$$ca = b \tag{1.2}$$

□

Aleks_Kleyn@MailAPS.org.
http://AleksKleyn.dyndns-home.com:4080/.
http://sites.google.com/site/AleksKleyn/.
http://arxiv.org/a/kleyn_a_1.
http://AleksKleyn.blogspot.com/.

It is evident that there is symmetry between definitions 1.1 and 1.2. The difference between left and right divisors is also evident since the product is noncommutative. However we can consider a definition generalizing definitions 1.1 and 1.2.

We may consider quotient of $b$ divided by $a$ as tuple of numbers $c$, $d$ such that[1]

$$cad = b$$

However we may consider division from another point of view.

Equations (1.1), (1.2) are examples of linear maps of $D$-algebra $A$. For commutative product, the equation (1.1) is only definition of linear map. For noncommutative product, $A \otimes A$-number generates linear map of $D$-algebra $A$. For instance, the tensor $c \otimes d$ is $A \otimes A$-number. Thus, we get the following definition.[2]

**Definition 1.3.** $A$-number $a$ is divisor of $A$-number $b$, if there exists $A \otimes A$-number $c$ such that

$$c \circ a = b$$

$A \otimes A$-number $c$ is called **quotient** of $A$-number $b$ divided by $A$-number $a$. □

In the paper [2], I considered equation

$$c \circ x = b$$

In this paper, I consider equation

$$x \circ a = b$$

$D$-algebra $A$ is called **division algebra**, if for any $A$-number $a \neq 0$ there exists $A$-number $a^{-1}$.

## 2. Conventions

**Convention 2.1.** *Element of $D$-algebra $A$ is called $A$-number. For instance, complex number is also called $C$-number, and quaternion is called $H$-number.* □

**Convention 2.2.** *Let $A$ be $\Omega_1$-algebra. Let $B$ be $\Omega_2$-algebra. Notation*

$$A \longrightarrow_* B$$

*means that there is representation of $\Omega_1$-algebra $A$ in $\Omega_2$-algebra $B$.* □

---

[1]Here somebody can argue that division is the inverse operation of multiplication. If product is commutative, then it follows from expression

$$ca = b$$

that
- $c$ is factor; the product of $c$ and $a$ is equal to $b$.
- $c$ is quotient of $b$ divided by $a$.

Thus, for commutative product, definitions of factor and quotient coincide. For noncommutative product, we distinguish right and left divisors; therefore, we distinguish left and right quotient.

[2]Since the map

(1.3) $$a \to cad$$

is bilinear, then, according to the theorem [3]-3.6.4, we may consider the map

$$a \to (c \otimes d) \circ a$$

instead of the map (1.3)

**Convention 2.3.** *Let A be associative D-algebra. The representation*

$$A \otimes A \xrightarrow{f}_* A \quad f(p) : a \to p \circ a$$

*of D-module $A \otimes A$ generates the set of linear maps. This representation generates product $\circ$ in D-module $A \otimes A$ according to rule*

$$(p \circ q) \circ a = p \circ (q \circ a)$$

□

Without a doubt, the reader may have questions, comments, objections. I will appreciate any response.

## 3. Geometry of Quotients

**Theorem 3.1.** *Let real field R be subfield of the center $Z(D)$ of the ring D. Then, for any A-numbers a, b, the set of quotients is convex.*

*Proof.* Let $A \otimes A$-number $c$ be quotient of $A$-number $b$ divided by $A$-number $a$

(3.1) $$c \circ a = b$$

Let $A \otimes A$-number $d$ be quotient of $A$-number $b$ divided by $A$-number $a$

(3.2) $$d \circ a = b$$

Then, for any $t \in R$, $0 \leq t \leq 1$, from (3.1), (3.2), it follows that

$$(tc + (1-t)d) \circ a = (tc) \circ a + ((1-t)d) \circ a$$
$$= t(c \circ a) + (1-t)(d \circ a) = tb + (1-t)b = b$$

Therefore, $A \otimes A$-number $tc + (1-t)d$ is quotient of $A$-number $b$ divided by $A$-number $a$. □

**Theorem 3.2.** *Let $A \otimes A$-numbers c, d be quotients of A-number b divided by A-number a*

(3.3) $$c \circ a = b$$

(3.4) $$d \circ a = b$$

*Then $A \otimes A$-number $c - d$ is quotient of A-number 0 divided by A-number a*

(3.5) $$(c - d) \circ a = 0$$

*Proof.* From (3.3), (3.4), it follows that

(3.6) $$(c - d) \circ a = c \circ a - d \circ a = b - b = 0$$

From (3.6), it follows that $A$-number $c - d$ is quotient of $A$-number 0 divided by $A$-number $a$. □

*Remark* 3.3. This point in the paper is interesting since it is easy to make mistake here and this mistake can be proved as a theorem. After we proved theorems 3.1, 3.2, the following question arises. If the difference of two quotients of $A$-number $b$ divided by $A$-number $a$ is quotient of $A$-number $0$ divided by $A$-number $a$, then what is the structure of the set of quotient $A$-number $0$ divided by any $A$-number? Or, equivalently, what is the structure of the set of quotients of $A$-number $b$ divided by $A$-number $a$?

Initially I wanted to prove the following statement.

*Let A-number a be neither a left nor right zero divisor. Then the quotient of A-number* 0 *divided by A-number a has either the form* 0*ac, or the form ca*0.

The proof is very simple. Since $c \ne 0$ in the expression *cad*, then $d = 0$.

When I decided to write down the quotient as

$$c \otimes 0 + 0 \otimes d$$

I asked myself whether I was right. Indeed, according to theorems [1]-13.2.2, [5]-9.2.2, there exists nontrivial representation of zero tensor. □

## 4. Division in $D$-Algebra

**Theorem 4.1.** *Let $D$-algebra $A$ be division algebra. Then, for any A-numbers $a$, $b$, there exist A-numbers $c$, $d$ such that*

$$cad = b$$

*Proof.* To prove the theorem it suffices to put

$$c = ba^{-1} \quad d = e$$

or

$$c = e \quad d = a^{-1}b$$

□

**Theorem 4.2.** *Let A-number $a$ divide A-number $b$. Let A-number $b$ divide A-number $c$. Then A-number $a$ divides A-number $c$.*

*Proof.* According to the definition 1.3, since A-number $a$ divides A-number $b$, then there exists $A \otimes A$-number $p$ such that

(4.1) $$b = p \circ a$$

According to the definition 1.3, since A-number $b$ divides A-number $c$, then there exists $A \otimes A$-number $q$ such that

(4.2) $$c = q \circ b$$

From equations (4.1), (4.2), it follows that

(4.3) $$c = q \circ (p \circ a) = (q \circ p) \circ a$$

According to the definition 1.3, from the equation (4.3), it follows that A-number $a$ divides A-number $c$. □

**Theorem 4.3.** *Let A-number $a$ divide A-number $b$. Let A-number $a$ divide A-number $c$. Then A-number $a$ divides A-number $b + c$.*

*Proof.* According to the definition 1.3, since A-number $a$ divides A-number $b$, then there exists $A \otimes A$-number $p$ such that

(4.4) $$b = p \circ a$$

According to the definition 1.3, since A-number $a$ divides A-number $c$, then there exists $A \otimes A$-number $q$ such that

(4.5) $$c = q \circ a$$

From equations (4.4), (4.5), it follows that

(4.6) $$b + c = p \circ a + q \circ a = (p + q) \circ a$$

According to the definition 1.3, from the equation (4.6), it follows that $A$-number $a$ divides $A$-number $b + c$. □

**Theorem 4.4.** *Let $A$-number $a$ divide $A$-number $b$. Let $c$ be $A \otimes A$-number. Then $A$-number $a$ divides $A$-number $c \circ b$.*

*Proof.* According to the definition 1.3, since $A$-number $a$ divides $A$-number $b$, then there exists $A \otimes A$-number $p$ such that

$$b = p \circ a \tag{4.7}$$

From the equation (4.7), it follows that

$$c \circ b = c \circ (p \circ a) = (c \circ p) \circ a \tag{4.8}$$

According to the definition 1.3, from the equation (4.8), it follows that $A$-number $a$ divides $A$-number $c \circ b$. □

## 5. Division with Remainder

**Definition 5.1.** *Let division in the $D$-algebra $A$ is not always defined. **$A$-number $a$ divides $A$-number $b$ with remainder**, if the following equation is true*

$$c \circ a + f = b \tag{5.1}$$

*$A \otimes A$-number $c$ is called **quotient** of $A$-number $b$ divided by $A$-number $a$. $A$-number $f$ is called **remainder of the division** of $A$-number $b$ by $A$-number $a$.* □

**Theorem 5.2.** *For any $A$-numbers $a$, $b$, there exist quotient and remainder of the division of $A$-number $b$ by $A$-number $a$.*

*Proof.* To prove the theorem, it is enough to assume that

$$c = a \otimes a$$
$$f = b - (a \otimes a) \circ a = b - aaa$$

□

From the theorem 5.2, it follows that, for given $A$-numbers $a$, $b$, a representation (5.1) is not unique. If we divide with remainder $A$-number $f$ by $A$-number $a$, then we get

$$c' \circ a + f' = f \tag{5.2}$$

From (5.1), (5.2), it follows that

$$c \circ a + c' \circ a + f' = (c + c') \circ a + f' = b \tag{5.3}$$

How can we compare representations (5.1), (5.3), when order is not defined in $D$-algebra $A$?

Therefore, the set of remainders of the division of $A$-number $b$ by $A$-number $a$ has form

$$A\{b, a\} = \{f \in A : f = b - c \circ a, c \in A \otimes A\} \tag{5.4}$$

If $0 \in A\{b, a\}$, then choice of remainder and corresponding quotient is evident. In general, a choice of a representative of the set $A\{b, a\}$ depends on properties of the algebra. In the algebra $N$ of integers, we choose the smallest positive number. In the algebra $A[x]$ of polynomials, we choose polynomial wich has the power less then power of divisor. If there is norm in algebra $A$, (for instance, the algebra of integer quaternions), then we can chose $A$-number with the smallest norm as remainder.

**Theorem 5.3.** *Since*

$$(5.5) \qquad A\{b,a\} \cap A\{c,a\} \neq \emptyset$$

*then*

  5.3.1: *$A$-number $a$* **divides $A$-number $b-c$ without remainder**.
  5.3.2: $A\{b,a\} = A\{c,a\}$

*Proof.* From the statement (5.5), it follows that there exists $A$-number

$$(5.6) \qquad d \in A\{b,a\} \cap A\{c,a\}$$

From the statement (5.6) and from the definition (5.4), it follows that there exist $A \otimes A$-numbers $f$, $g$ such that

$$(5.7) \qquad b = f \circ a + d$$

$$(5.8) \qquad c = g \circ a + d$$

From equations (5.7), (5.8), it follows that

$$(5.9) \qquad b - c = f \circ a - g \circ a = (f - g) \circ a$$

The statement 5.3.1 follows from the equation (5.9).

Let $m \in A\{b,a\}$. From the definition (5.4), it follows that there exist $A \otimes A$-number $n$ such that

$$(5.10) \qquad b = n \circ a + m$$

From equations (5.7), (5.10), it follows that

$$(5.11) \qquad f \circ a + d = n \circ a + m$$

From the equation (5.11), it follows that

$$(5.12) \qquad d = n \circ a - f \circ a + m = (n - f) \circ a + m$$

From equations (5.8), (5.12), it follows that

$$(5.13) \qquad c = g \circ a + (n - f) \circ a + m = (g + n - f) \circ a + m$$

From the definition (5.4) and the equation (5.13), it follows that there exist $m \in A\{c,a\}$. Therefore,

$$(5.14) \qquad A\{b,a\} \subseteq A\{c,a\}$$

The same way, we prove the statement

$$(5.15) \qquad A\{c,a\} \subseteq A\{b,a\}$$

The statement 5.3.2 folows from statements (5.14), (5.15). □

From the theorem 5.3, it follows that, for given $A$-number $a$, the family of sets $A\{b,a\}$ generates equivalence mod $a$.

**Definition 5.4.** *We define* **canonical remainder** $b \bmod a$ **of the division** *of $A$-number $b$ by $A$-number $a$ as selected element of the set $A\{b,a\}$. The representation*

$$c \circ a + (b \bmod a) = b$$

*of division with remainder is called* **canonical**. □

At first glance, the choice of $A$-number $b \bmod a$ is arbitrary. However we can define the natural constrains of the arbitrary choice.

**Theorem 5.5.** *If we define sum on the set $A/\bmod a$ according to the rule*

(5.16) $$b \bmod a + c \bmod a = (b+c) \bmod a$$

*then the set $A/\bmod a$ is Abelian group.*

*Proof.* According to the definition 5.4, there exist $A \otimes A$-numbers $p$, $q$ such that

(5.17) $$\begin{aligned} b &= p \circ a + b \bmod a \\ c &= q \circ a + c \bmod a \end{aligned}$$

From (5.17), it follows that

(5.18) $$\begin{aligned} b + c &= p \circ a + b \bmod a + q \circ a + c \bmod a \\ &= (p+q) \circ a + b \bmod a + c \bmod a \end{aligned}$$

From equations (5.4), (5.18), it follows that

(5.19) $$b \bmod a + c \bmod a \in A\{b+c, a\}$$

From the statement (5.19) and from the definition 5.4, it follows that sum (5.16) is well defined.

We verify commutativity of the sum (5.16) directly. □

**Theorem 5.6.** *The representation*

$$D \dashrightarrow A/\bmod a$$

*of ring $D$ in Abelian group $A/\bmod a$ defined by the equation*

(5.20) $$d(b \bmod a) = (db) \bmod a$$

*generates $D$-module $A/\bmod a$.*

*Proof.* According to the definition 5.4, there exist $A \otimes A$-number $p$ such that

(5.21) $$b = p \circ a + b \bmod a$$

From (5.21), it follows that

(5.22) $$db = d(p \circ a) + d(b \bmod a) = (dp) \circ a + d(b \bmod a)$$

From equations (5.4), (5.22), it follows that

(5.23) $$d(b \bmod a) \in A\{db, a\}$$

From the statement (5.23) and from the definition 5.4, it follows that representation (5.20) is well defined. □

**Theorem 5.7.** *If we define product in $D$-module $A/\bmod a$ according to the rule*

(5.24) $$(b \bmod a)(c \bmod a) = (bc) \bmod a$$

*then $D$-module $A/\bmod a$ is $D$-algebra.*

*Proof.* According to the definition 5.4, there exist $A \otimes A$-numbers $p$, $q$ such that

$$b = p \circ a + b \bmod a \tag{5.25}$$
$$c = q \circ a + c \bmod a$$

From (5.25), it follows that

$$\begin{aligned}
bc &= (p \circ a + b \bmod a)(q \circ a + c \bmod a) \\
&= (p \circ a)(q \circ a + c \bmod a) + (b \bmod a)(q \circ a + c \bmod a) \\
&= (p \circ a)b + (b \bmod a)(q \circ a) + (b \bmod a)(c \bmod a) \\
&= ((1 \otimes b) \circ p) \circ a + (((b \bmod a) \otimes 1) \circ q) \circ a + (b \bmod a)(c \bmod a) \\
&= ((1 \otimes b) \circ p + ((b \bmod a) \otimes 1) \circ q) \circ a + (b \bmod a)(c \bmod a)
\end{aligned} \tag{5.26}$$

From equations (5.4), (5.26), it follows that

$$(b \bmod a)(c \bmod a) \in A\{bc, a\} \tag{5.27}$$

From the statement (5.27) and from the definition 5.4, it follows that product (5.24) is well defined. □

**Theorem 5.8.** *Let*

$$p \circ b + q = a \tag{5.28}$$

*be canonical representation of division with remainder of $A$-number $a$ by $A$-number $b$. Let*

$$t \circ c + s = b \tag{5.29}$$

*be canonical representation of division with remainder of $A$-number $b$ by $A$-number $c$.*

  5.8.1: *Let*

$$u \circ c + v = p \circ s + q \tag{5.30}$$

  *be canonical representation of division with remainder of $A$-number $p \circ s + q$ by $A$-number $c$.*

  5.8.2: *Then the canonical representation of division with remainder of $A$-number $a$ by $A$-number $c$ has form*

$$(p \circ t + u) \circ c + v = a \tag{5.31}$$

*Proof.* From equations (5.28), (5.29), it follows that

$$a = p \circ (t \circ c + s) + q = p \circ t \circ c + p \circ s + q \tag{5.32}$$

The equation (5.32) is representation of division with remainder of $A$-number $a$ by $A$-number $c$. From equations (5.30), (5.32), it follows that

$$a = p \circ t \circ c + u \circ c + v \tag{5.33}$$

The equation (5.31) follows from the equation (5.33). From the statement 5.8.1 and from the definition 5.4, it follows that

$$v = (p \circ s + q) \bmod c \tag{5.34}$$

From equations (5.4), (5.31), it follows that

$$v \in A\{a, c\} \tag{5.35}$$

From the statement (5.35) and from the definition 5.4, it follows that

(5.36) $$v = a \bmod c$$

The statement 5.8.2 follows from the statement (5.36) and from the definition 5.4. □

## 6. Highest Common Factor

**Definition 6.1.** *$A$-number $c$ is called* **common factor** *of $A$-numbers $a$ and $b$, if $A$-number $c$ divides each of $A$-numbers $a$ and $b$. If $A$-numbers $a$ and $b$ are not unit divisors and any common factor of $A$-numbers $a$ and $b$ is not unit divisor, then $A$-numbers $a$ and $b$ are called* **relatively prime**. □

**Definition 6.2.** *$A$-number $c$ is called* **highest common factor** *of $A$-numbers $a$ and $b$, if $A$-number $c$ is common factor of $A$-numbers $a$ and $b$ and any common factor $d$ of $A$-numbers $a$ and $b$ divides $A$-number $c$.* □

## 7. Prime $A$-number

**Definition 7.1.** *Let $A$-number $b$ be not unit divisor of $D$-algebra $A$. $A$-number $b$ is called* **prime**, *if any divisor $a$ of $A$-number $b$ satisfies one of the following conditions.*

7.1.1: *$A$-number $a$ is unit divisor.*
7.1.2: *Quotient of $A$-number $b$ divided by $A$-number $a$ is unit divisor of $D$-algebra. $A \otimes A$.*

□

**Theorem 7.2.** *Let $D$-algebra $A$ is division algebra. Let $A[x]$ be algebra of polynomials over $D$-algebra $A$. Polynomial of power $1$ is prime $A[x]$-number.*

*Proof.* Let $p$ be polynomial of power 1. Let $q$ be polynomial. Let

(7.1) $$p = (r_1 \otimes r_2) \circ q = r_1 q r_2$$

be canonical representation of division with remainder of polynomial $p$ over polynomial $q$. Here $r_1$, $r_2$ are polynomials. According to the theorems [4]-5.9, [6]-20,

(7.2) $$\deg r_1 + \deg q + \deg r_2 = \deg p = 1$$

From the equation (7.2), it follows that only one polynomial $q$, $r_1$, $r_2$ has power 1, and other two polynomials are $A$-numbers.

- Since the polynomial $q$ is $A$-number, then the polynomial $q$ is unit divisor and satisfies the condition 7.1.1.
- Since $q$ is polynomial of power 1, then, according to the theorems [6]-28, [4]-6.10,

(7.3) $$p = r \circ q$$

  where $r$ is $A \otimes A$-number. Since $p$ and $q$ are polynomials of power 1, then, according to the theorems [6]-28, [4]-6.10,

(7.4) $$q = r' \circ p + s$$

  where $r'$ is $A \otimes A$-number and $s$ is $A$-number. From equations (7.3), (7.4), it follows that

(7.5) $$p = r \circ r' \circ p + r \circ s$$

...

From the equation (7.5), it follows that

(7.6) $$r \circ r' = 1 \otimes 1 \quad r \circ s = 0$$

From the equation (7.6), it follows that $A \otimes A$-number $r$ is unit divisor and $s = 0$. Therefore, the polynomial $q$ satisfies the condition 7.1.2.

According to the definition 7.1, the polynomial $p$ is prime $A[x]$-number. $\square$

## 8. References

[1] Aleks Kleyn, Lectures on Linear Algebra over Division Ring,
eprint arXiv:math.GM/0701238 (2010)

[2] Aleks Kleyn, Linear Equation in Finite Dimensional Algebra,
eprint arXiv:0912.4061 (2010)

[3] Aleks Kleyn, Linear Maps of Free Algebra,
eprint arXiv:1003.1544 (2010)

[4] Aleks Kleyn, Polynomial over Associative $D$-Algebra,
eprint arXiv:1302.7204 (2013)

[5] Aleks Kleyn.
Linear Algebra over Division Ring: Vector Space.
CreateSpace, 2014; ISBN-13: 978-1499324006

[6] Aleks Kleyn, Polynomial over Associative $D$-Algebra.
Clifford Analysis, Clifford Algebras and their applications, Vol 2, Issue 2,
pages 97 - 115, 2013

# 9. Index

## 10. Special Symbols and Notations

# Деление в ассоциативной $D$-алгебре

Александр Клейн

**Аннотация.** Из симметрии между определениями левого и правого делителей в ассоциативной $D$-алгебре $A$ следует возможность определить частное как $A \otimes A$-число. В статье рассмотрены деление и деление с остатком. Я рассмотрел также понятие простого $A$-числа.

Содержание

## 1. Предисловие

Пусть $D$ - коммутативное кольцо. Мы будем предполагать, что $D$-алгебра $A$ ассоциативна.

**Определение 1.1.** $A$-число $a$ называется левым делителем $A$-числа $b$, если существует $A$-число $c$ такое, что

$$ac = b \quad (1.1)$$

□

**Определение 1.2.** $A$-число $a$ называется правым делителем $A$-числа $b$, если существует $A$-число $c$ такое, что

$$ca = b \quad (1.2)$$

Aleks_Kleyn@MailAPS.org.
http://AleksKleyn.dyndns-home.com:4080/.
http://sites.google.com/site/AleksKleyn/.
http://arxiv.org/a/kleyn_a_1.
http://AleksKleyn.blogspot.com/.

□

Симметрия между определениями 1.1 и 1.2 очевидна. Также как очевидно различие между левым и правым делителями в связи с некоммутативностью произведения. Однако мы можем рассмотреть определение, обобщающее определения 1.1 и 1.2.

Мы можем рассматривать частное от деления $b$ на $a$ как пару чисел $c$, $d$ таких, что[1]

$$cad = b$$

Однако мы можем рассмотреть операцию деления с другой точки зрения.

Равенства (1.1), (1.2) являются примерами линейных отображений $D$-алгебры $A$. В коммутативном случае равенство (1.1) является единственным определением линейного отображения. В некоммутативном случае, линейное отображение $D$-алгебры $A$ порождается $A \otimes A$-числом. Например, тензор $c \otimes d$ является $A \otimes A$-числом. Таким образом, мы получаем следующее определение.[2]

**Определение 1.3.** $A$-число $a$ называется делителем $A$-числа $b$, если существует $A \otimes A$-число $c$ такое, что

$$c \circ a = b$$

$A \otimes A$-число $c$ называется **частным от деления** $A$-числа $b$ на $A$-число $a$.  □

В статье [2] я рассмотрел уравнение

$$c \circ x = b$$

В этой статье я рассматриваю уравнение

$$x \circ a = b$$

$D$-алгебра $A$ называется **алгеброй с делением**, если для любого $A$-числа $a \neq 0$ существует $A$-число $a^{-1}$.

## 2. Соглашения

**Соглашение 2.1.** *Элемент $D$-алгебры $A$ называется $A$-**числом**. Например, комплексное число также называется $C$-числом, а кватернион называется $H$-числом.*  □

---

[1]Здесь можно возразить, что деление - это операция, обратная умножению. В коммутативном случае, из выражения

$$ca = b$$

следует, что
- $c$ является множителем, произведение $c$ и $a$ равно $b$.
- $c$ является частным деления $b$ на $a$.

Таким образом, в коммутативном случае определения множителя и частного совпадают. В некоммутативном случае, мы различаем правый и левый делители; следовательно, мы различаем левое и правое частное.

[2]Так как отображение

(1.3) $$a \to cad$$

билинейно, то, согласно теореме [3]-3.6.4, мы можем рассматривать отображение

$$a \to (c \otimes d) \circ a$$

вместо отображения (1.3).

**Соглашение 2.2.** *Пусть $A$ - $\Omega_1$-алгебра. Пусть $B$ - $\Omega_2$-алгебра. Запись*

$$A \xrightarrow{\ \ *\ } B$$

*означает, что определено представление $\Omega_1$-алгебры $A$ в $\Omega_2$-алгебре $B$.* □

**Соглашение 2.3.** *Пусть $A$ ассоциативная $D$-алгебра. Представление*

$$A \otimes A \xrightarrow{\ \ f\ }_{*} A \quad f(p) : a \to p \circ a$$

*$D$-модуля $A \otimes A$ порождает множество линейных отображений. Это представление порождает произведение $\circ$ в $D$-модуле $A \otimes A$ согласно правилу*

$$(p \circ q) \circ a = p \circ (q \circ a)$$

□

Без сомнения, у читателя могут быть вопросы, замечания, возражения. Я буду признателен любому отзыву.

## 3. Геометрия частных

**Теорема 3.1.** *Пусть поле действительных чисел $R$ является подполем центра $Z(D)$ кольца $D$. Тогда для любых $A$-чисел $a$, $b$ множество частных выпукло.*

*Доказательство.* Пусть $A \otimes A$-число $c$ является частным от деления $A$-числа $b$ на $A$-число $a$

$$(3.1) \qquad c \circ a = b$$

Пусть $A \otimes A$-число $d$ является частным от деления $A$-числа $b$ на $A$-число $a$

$$(3.2) \qquad d \circ a = b$$

Тогда для любого $t \in R$, $0 \le t \le 1$, из (3.1), (3.2) следует, что

$$(tc + (1-t)d) \circ a = (tc) \circ a + ((1-t)d) \circ a$$
$$= t(c \circ a) + (1-t)(d \circ a) = tb + (1-t)b = b$$

Следовательно, $A \otimes A$-число $tc + (1-t)d$ является частным от деления $A$-числа $b$ на $A$-число $a$. □

**Теорема 3.2.** *Пусть $A \otimes A$-числа $c$, $d$ являются частными от деления $A$-числа $b$ на $A$-число $a$*

$$(3.3) \qquad c \circ a = b$$

$$(3.4) \qquad d \circ a = b$$

*Тогда $A \otimes A$-число $c - d$ является частным от деления $A$-числа $0$ на $A$-число $a$*

$$(3.5) \qquad (c - d) \circ a = 0$$

*Доказательство.* Из (3.3), (3.4) следует, что

$$(3.6) \qquad (c - d) \circ a = c \circ a - d \circ a = b - b = 0$$

Из (3.6) следует, что $A$-число $c - d$ является частным от деления $A$-числа $0$ на $A$-число $a$. □

*Замечание* 3.3. Это место в статье очень интересно тем, что здесь очень легко совершить ошибку, которую легко доказать как теорему. После доказательства теорем 3.1, 3.2, естественно возникает следующий вопрос. Если разность двух частных от деления $A$-числа $b$ на $A$-число $a$ является частным от деления $A$-числа $0$ на $A$-число $a$, то какова структура множества частных от деления $A$-числа $0$? Или, что тоже самое, какова структура множества частных от деления $A$-числа $b$ на $A$-число $a$?

Вначале я хотел доказать следующее утверждение.

*Пусть $A$-число $a$ не является ни левым, ни правым делителем нуля. Тогда частное от деления $A$-числа $0$ на $A$-число $a$ имеет либо вид $0ac$, либо вид $ca0$.*

Доказательство очень просто. Если $c \ne 0$ в выражении $cad$, то $d = 0$.

Когда я решил записать частное в виде

$$c \otimes 0 + 0 \otimes d$$

я подумал прав ли я. Действительно, согласно теоремам [1]-13.2.2, [5]-9.2.2, существует нетривиальная запись нулевого тензора. □

## 4. Деление в $D$-алгебре

**Теорема 4.1.** *Если $D$-алгебра $A$ является алгеброй с делением, то для любых $A$-чисел $a$, $b$ существуют $A$-числа $c$, $d$ такие, что*

$$cad = b$$

*Доказательство.* Для доказательства теоремы достаточно положить

$$c = ba^{-1} \quad d = e$$

или

$$c = e \quad d = a^{-1}b$$

□

**Теорема 4.2.** *Пусть $A$-число $a$ делит $A$-число $b$. Пусть $A$-число $b$ делит $A$-число $c$. Тогда $A$-число $a$ делит $A$-число $c$.*

*Доказательство.* Согласно определению 1.3, так как $A$-число $a$ делит $A$-число $b$, то существует $A \otimes A$-число $p$ такое, что

(4.1) $$b = p \circ a$$

Согласно определению 1.3, так как $A$-число $b$ делит $A$-число $c$, то существует $A \otimes A$-число $q$ такое, что

(4.2) $$c = q \circ b$$

Из равенств (4.1), (4.2) следует, что

(4.3) $$c = q \circ (p \circ a) = (q \circ p) \circ a$$

Согласно определению 1.3, из равенства (4.3) следует, что $A$-число $a$ делит $A$-число $c$. □

**Теорема 4.3.** *Пусть $A$-число $a$ делит $A$-число $b$. Пусть $A$-число $a$ делит $A$-число $c$. Тогда $A$-число $a$ делит $A$-число $b + c$.*

*Доказательство.* Согласно определению 1.3, так как $A$-число $a$ делит $A$-число $b$, то существует $A \otimes A$-число $p$ такое, что

$$b = p \circ a \tag{4.4}$$

Согласно определению 1.3, так как $A$-число $a$ делит $A$-число $c$, то существует $A \otimes A$-число $q$ такое, что

$$c = q \circ a \tag{4.5}$$

Из равенств (4.4), (4.5) следует, что

$$b + c = p \circ a + q \circ a = (p + q) \circ a \tag{4.6}$$

Согласно определению 1.3, из равенства (4.6) следует, что $A$-число $a$ делит $A$-число $b + c$. $\square$

**Теорема 4.4.** *Пусть $A$-число $a$ делит $A$-число $b$. Пусть $c$ является $A \otimes A$-числом. Тогда $A$-число $a$ делит $A$-число $c \circ b$.*

*Доказательство.* Согласно определению 1.3, так как $A$-число $a$ делит $A$-число $b$, то существует $A \otimes A$-число $p$ такое, что

$$b = p \circ a \tag{4.7}$$

Из равенств (4.7) следует, что

$$c \circ b = c \circ (p \circ a) = (c \circ p) \circ a \tag{4.8}$$

Согласно определению 1.3, из равенства (4.8) следует, что $A$-число $a$ делит $A$-число $c \circ b$. $\square$

## 5. Деление с остатком

**Определение 5.1.** *Пусть деление в $D$-алгебре $A$ не всегда определено. **$A$-число $a$ делит $A$-число $b$ с остатком**, если следующее равенство верно*

$$c \circ a + f = b \tag{5.1}$$

*$A \otimes A$-число $c$ называется **частным от деления** $A$-числа $b$ на $A$-число $a$. $A$-число $d$ называется **остаток от деления** $A$-числа $b$ на $A$-число $a$.* $\square$

**Теорема 5.2.** *Для любых $A$-чисел $a$, $b$, существует частное и остаток от деления $A$-числа $b$ на $A$-число $a$.*

*Доказательство.* Для доказательства теоремы достаточно положить

$$c = a \otimes a$$
$$f = b - (a \otimes a) \circ a = b - aaa$$

$\square$

Из теоремы 5.2 следует, что для заданных $A$-чисел $a$, $b$ представление (5.1) определено не однозначно. Если мы поделим с остатком $A$-число $f$ на $A$-число $a$, то мы получим

$$c' \circ a + f' = f \tag{5.2}$$

Из (5.1), (5.2) следует, что

$$c \circ a + c' \circ a + f' = (c + c') \circ a + f' = b \tag{5.3}$$

Как мы можем сравнить представления (5.1), (5.3), если отношение порядка не определено в $D$-алгебре $A$?

Следовательно, множество остатков от деления $A$-числа $b$ на $A$-число $a$ имеет вид

$$(5.4) \qquad A\{b,a\} = \{f \in A : f = b - c \circ a, c \in A \otimes A\}$$

Если $0 \in A\{b,a\}$, то выбор остатка и соответствующего частного очевиден. Вообще говоря, выбор представителя множества $A\{b,a\}$ зависит от свойств алгебры. В алгебре $N$ целых чисел мы выбираем наименьшее положительное число. В алгебре $A[x]$ многочленов мы выбираем многочлен степени меньше степени делителя. Если в алгебре $A$ определена норма (например, алгебра целых кватернионов), то в качестве остатка мы можем выбрать $A$-число с наименьшей нормой.

**Теорема 5.3.** *Если*

$$(5.5) \qquad A\{b,a\} \cap A\{c,a\} \neq \emptyset$$

*то*

  5.3.1: *$A$-число $a$* **делит $A$-число $b - c$ без остатка**.
  5.3.2: $A\{b,a\} = A\{c,a\}$

*Доказательство.* Из утверждения (5.5) следует, что существует $A$-число

$$(5.6) \qquad d \in A\{b,a\} \cap A\{c,a\}$$

Из утверждения (5.6) и определения (5.4) следует, что существуют $A \otimes A$-числа $f$, $g$ такие, что

$$(5.7) \qquad b = f \circ a + d$$

$$(5.8) \qquad c = g \circ a + d$$

Из равенств (5.7), (5.8) следует, что

$$(5.9) \qquad b - c = f \circ a - g \circ a = (f - g) \circ a$$

Утверждение 5.3.1 является следствием равенства (5.9).

Пусть $m \in A\{b,a\}$. Из определения (5.4) следует, что существует $A \otimes A$-число $n$ такие, что

$$(5.10) \qquad b = n \circ a + m$$

Из равенств (5.7), (5.10) следует, что

$$(5.11) \qquad f \circ a + d = n \circ a + m$$

Из равенства (5.11) следует, что

$$(5.12) \qquad d = n \circ a - f \circ a + m = (n - f) \circ a + m$$

Из равенств (5.8), (5.12) следует, что

$$(5.13) \qquad c = g \circ a + (n - f) \circ a + m = (g + n - f) \circ a + m$$

Из определения (5.4) и равенства (5.13) следует, что существует $m \in A\{c,a\}$. Следовательно,

$$(5.14) \qquad A\{b,a\} \subseteq A\{c,a\}$$

Аналогичным образом мы доказываем утверждение

(5.15) $$A\{c, a\} \subseteq A\{b, a\}$$

Утверждение 5.3.2 является следствием утверждений (5.14), (5.15). □

Из теоремы 5.3 следует, что для данного $A$-числа $a$ семейство множеств $A\{b, a\}$ порождает отношение эквивалентности $\mod a$.

**Определение 5.4.** *Определим* **каноническое частное** $b \mod a$ *от деления $A$-числа $b$ на $A$-число $a$ как выбранный элемент множества $A\{b, a\}$. Представление*

$$c \circ a + (b \mod a) = b$$

*деления с остатком называется* **каноническим**. □

На первый взгляд, выбор $A$-числа $b \mod a$ произволен. Однако мы можем определить естественные границы этого произвола.

**Теорема 5.5.** *Если мы определим сложение на множестве $A/\mod a$ согласно правилу*

(5.16) $$b \mod a + c \mod a = (b + c) \mod a$$

*то множество $A/\mod a$ является абелевой группой.*

*Доказательство.* Согласно определению 5.4, существуют $A \otimes A$-числа $p$, $q$ такие, что

(5.17) $$\begin{aligned} b &= p \circ a + b \mod a \\ c &= q \circ a + c \mod a \end{aligned}$$

Из (5.17) следует, что

(5.18) $$\begin{aligned} b + c &= p \circ a + b \mod a + q \circ a + c \mod a \\ &= (p + q) \circ a + b \mod a + c \mod a \end{aligned}$$

Из равенств (5.4), (5.18) следует, что

(5.19) $$b \mod a + c \mod a \in A\{b + c, a\}$$

Из утверждения (5.19) и определения 5.4 следует корректность определения (5.16) суммы.

Коммутативность суммы (5.16) доказывается непосредственной проверкой. □

**Теорема 5.6.** *Представление*

$$D \xrightarrow{\ *\ } A/\mod a$$

*кольца $D$ в абелевой группе $A/\mod a$ определённое равенством*

(5.20) $$d(b \mod a) = (db) \mod a$$

*порождает $D$-модуль $A/\mod a$.*

*Доказательство.* Согласно определению 5.4, существует $A \otimes A$-число $p$ такое, что

(5.21) $$b = p \circ a + b \bmod a$$

Из (5.21) следует, что

(5.22) $$db = d(p \circ a) + d(b \bmod a) = (dp) \circ a + d(b \bmod a)$$

Из равенств (5.4), (5.22) следует, что

(5.23) $$d(b \bmod a) \in A\{db, a\}$$

Из утверждения (5.23) и определения 5.4 следует корректность определения (5.20) представления. □

**Теорема 5.7.** *Если мы определим умножение в $D$-модуле $A/\bmod a$ согласно правилу*

(5.24) $$(b \bmod a)(c \bmod a) = (bc) \bmod a$$

*то $D$-модуль $A/\bmod a$ является $D$-алгеброй.*

*Доказательство.* Согласно определению 5.4, существуют $A \otimes A$-числа $p$, $q$ такие, что

(5.25) $$\begin{aligned} b &= p \circ a + b \bmod a \\ c &= q \circ a + c \bmod a \end{aligned}$$

Из (5.25) следует, что

(5.26) $$\begin{aligned} bc &= (p \circ a + b \bmod a)(q \circ a + c \bmod a) \\ &= (p \circ a)(q \circ a + c \bmod a) + (b \bmod a)(q \circ a + c \bmod a) \\ &= (p \circ a)b + (b \bmod a)(q \circ a) + (b \bmod a)(c \bmod a) \\ &= ((1 \otimes b) \circ p) \circ a + (((b \bmod a) \otimes 1) \circ q) \circ a + (b \bmod a)(c \bmod a) \\ &= ((1 \otimes b) \circ p + ((b \bmod a) \otimes 1) \circ q) \circ a + (b \bmod a)(c \bmod a) \end{aligned}$$

Из равенств (5.4), (5.26) следует, что

(5.27) $$(b \bmod a)(c \bmod a) \in A\{bc, a\}$$

Из утверждения (5.27) и определения 5.4 следует корректность определения (5.24) произведения. □

**Теорема 5.8.** *Пусть*

(5.28) $$p \circ b + q = a$$

*является каноническим представлением деления с остатком $A$-числа $a$ на $A$-число $b$. Пусть*

(5.29) $$t \circ c + s = b$$

*является каноническим представлением деления с остатком $A$-числа $b$ на $A$-число $c$.*

5.8.1: *Пусть*

(5.30) $$u \circ c + v = p \circ s + q$$

*является каноническим представлением деления с остатком $A$-числа $p \circ s + q$ на $A$-число $c$.*

5.8.2: *Тогда каноническое представление деления с остатком $A$-числа $a$ на $A$-число $c$ имеет вид*

$$(p \circ t + u) \circ c + v = a \tag{5.31}$$

*Доказательство.* Из равенств (5.28), (5.29) следует, что

$$a = p \circ (t \circ c + s) + q = p \circ t \circ c + p \circ s + q \tag{5.32}$$

Равенство (5.32) является представлением деления с остатком $A$-числа $a$ на $A$-число $c$. Из равенств (5.30), (5.32) следует, что

$$a = p \circ t \circ c + u \circ c + v \tag{5.33}$$

Равенство (5.31) является следствием равенства (5.33). Из утверждения 5.8.1 и определения 5.4 следует, что

$$v = (p \circ s + q) \bmod c \tag{5.34}$$

Из равенств (5.4), (5.31) следует, что

$$v \in A\{a, c\} \tag{5.35}$$

Из утверждения (5.35) и определения 5.4 следует, что

$$v = a \bmod c \tag{5.36}$$

Утверждение 5.8.2 является следствием утверждения (5.36) и определения 5.4. $\square$

## 6. Наибольший общий делитель

**Определение 6.1.** *$A$-число $c$ называется **общим делителем** $A$-чисел $a$ и $b$, если $A$-число $c$ делит каждое из $A$-чисел $a$ и $b$. Если $A$-числа $a$ и $b$ не являются делителями единицы и любой общий делитель $A$-чисел $a$ и $b$ не является делителем единицы, то $A$-числа $a$ и $b$ называются **взаимно простыми**.* $\square$

**Определение 6.2.** *$A$-число $c$ называется **наибольшим общим делителем** $A$-чисел $a$ и $b$, если $A$-число $c$ является общим делителем $A$-чисел $a$ и $b$ и любой общий делитель $d$ $A$-чисел $a$ и $b$ делит $A$-число $c$.* $\square$

## 7. Простое $A$-число

**Определение 7.1.** *Пусть $A$-число $b$ не является делителем единицы $D$-алгебры $A$. $A$-число $b$ называется **простым**, если любой делитель $a$ $A$-числа $b$ удовлетворяет одному из следующих условий.*

7.1.1: *$A$-число $a$ является делителем единицы.*
7.1.2: *Частное от деления $A$-числа $b$ на $A$-число $a$ является делителем единицы $D$-алгебры $A \otimes A$.*

$\square$

**Теорема 7.2.** *Пусть $D$-алгебра $A$ является алгеброй с делением. Рассмотрим алгебру многочленов $A[x]$ над $D$-алгеброй $A$. Многочлен степени 1 является простым $A[x]$-числом.*

*Доказательство.* Пусть $p$ - многочлен степени 1. Пусть $q$ - многочлен. Пусть

$$(7.1) \qquad p = (r_1 \otimes r_2) \circ q = r_1 q r_2$$

каноническая форма деления с остатком многочлена $p$ на многочлен $q$. Здесь $r_1$, $r_2$ - многочлены. Согласно теоремам [4]-5.9, [6]-20,

$$(7.2) \qquad \deg r_1 + \deg q + \deg r_2 = \deg p = 1$$

Из равенства (7.2) следует, что только один многочлен $q$, $r_1$, $r_2$ имеет степень 1, и остальные два многочлена являются $A$-числами.

- Если многочлен $q$ является $A$-числом, то многочлен $q$ является делителем единицы и удовлетворяет условию 7.1.1.
- Если степень многочлена $q$ равна 1, то, согласно теоремам [4]-6.10, [6]-28,

  $$(7.3) \qquad p = r \circ q$$

  где $r$ является $A \otimes A$-числом. Так как $p$, $q$ - многочлены степени 1, то, согласно теоремам [4]-6.10, [6]-28,

  $$(7.4) \qquad q = r' \circ p + s$$

  где $r'$ является $A \otimes A$-числом и $s$ является $A$-числом. Из равенств (7.3), (7.4) следует, что

  $$(7.5) \qquad p = r \circ r' \circ p + r \circ s$$

  Из равенства (7.5) следует, что

  $$(7.6) \qquad r \circ r' = 1 \otimes 1 \quad r \circ s = 0$$

  Из равенства (7.6) следует, что $A \otimes A$-число $r$ является делителем единицы и $s = 0$. Следовательно, многочлен $q$ удовлетворяет условию 7.1.2.

Согласно определению 7.1, многочлен $p$ является простым $A[x]$-числом. $\square$

## 8. Список литературы

[1] Александр Клейн, Лекции по линейной алгебре над телом, eprint arXiv:math.GM/0701238 (2010)
[2] Александр Клейн, Линейное уравнение в конечномерной алгебре, eprint arXiv:0912.4061 (2010)
[3] Александр Клейн, Линейные отображения свободной алгебры, eprint arXiv:1003.1544 (2010)
[4] Александр Клейн, Многочлен над ассоциативной $D$-алгеброй, eprint arXiv:1302.7204 (2013)
[5] Александр Клейн. Линейная алгебра над телом: Векторное пространство. CreateSpace, 2014; ISBN-13: 978-1499323948
[6] Aleks Kleyn, Polynomial over Associative $D$-Algebra. Clifford Analysis, Clifford Algebras and their applications, Vol 2, Issue 2, pages 97 - 115, 2013

# 9. Предметный указатель

# 10. Специальные символы и обозначения

$b \bmod a$   каноническое частное от
деления  7